# $\mathcal{MID}$ AND SUBNORMAL SAFE QUOTIENTS FOR GEOMETRICALLY REGULAR WEIGHTED SHIFTS

CHAFIQ BENHIDA, RAÚL E. CURTO, AND GEORGE R. EXNER

Abstract. Geometrically regular weighted shifts (in short, GRWS) are those with weights $\alpha(N, D)$ given by $\alpha_n(N, D) = \sqrt{\frac{p^n + N}{p^n + D}}$, where $p > 1$ and $(N, D)$ is fixed in the open unit square $(-1, 1) \times (-1, 1)$. We study here the zone of pairs $(M, P)$ for which the weight $\frac{\alpha(N,D)}{\alpha(M,P)}$ gives rise to a moment infinitely divisible ($\mathcal{MID}$) or a subnormal weighted shift, and deduce immediately the analogous results for product weights $\alpha(N, D)\alpha(M, P)$, instead of quotients. Useful tools introduced for this study are a pair of partial orders on the GRWS.

## 1. Introduction and preliminaries

Let $\mathcal{H}$ be a separable, complex infinite dimensional Hilbert space and $\mathcal{L}(\mathcal{H})$ the algebra of bounded linear operators on $\mathcal{H}$. The current paper continues the study, initiated in [5], of certain weighted shifts in $\mathcal{L}(\mathcal{H})$ called the geometrically regular weighted shifts (henceforth, GRWS), with concentration on when their quotients and products belong to standard classes of interest, especially the $\mathcal{MID}$ and subnormal shifts (all definitions reviewed below). We introduce as well two useful partial orders on these weighted shifts which aid in the study, and recast some earlier results in terms of these orders.

Study of such quotients is pertinent because the resulting shifts add to the collections of reasonably tractable shifts which nonetheless extend the standard shift examples. Shifts are often used to begin the study of a new concept or definition, or as test operators for a new hypothesis. However, we have discovered that the most used ones – the Agler-type shifts related to the Bergman shift and the finitely atomic (recursively generated) shifts – are so special as to be seriously misleading. The GRWS themselves are such an addition to the example set; their quotients, in various classes of considerable interest as defined by a base point and other points in the parameter space, yield a "quotient safe zone" for some property, and they are a further such addition. Related to these considerations is the fact that $\mathcal{MID}$ shifts form a sophisticated yet tractable collection of operators, with many important properties, among others their stability under the Aluthge transform and Schur products [4]. As a result, any time two GRWS $W_\alpha$ and $W_\beta$ have a (Schur) $\mathcal{MID}$ quotient $W_{\frac{\alpha}{\beta}}$, all reasonable properties involving moment matrices or Agler sums enjoyed by $W_\beta$ will transfer to $W_\alpha$. This places the class of $\mathcal{MID}$ GRWS in a privileged position, and justifies our emphasis on safe quotients.

We first briefly recall the relevant definitions and set notation. For a more complete introduction and fuller discussion, see [5].

## 1.1. Unilateral weighted shifts and operator properties.

Set $\mathbb{N}_0 := \{0, 1, \ldots\}$ and let $\ell^2$ be the classical Hilbert space $\ell^2(\mathbb{N}_0)$ with canonical orthonormal basis $e_0, e_1, \ldots$. Let $\alpha : \alpha_0, \alpha_1, \ldots$ be a (bounded) non-negative **weight sequence** and $W_\alpha$ the weighted shift defined by $W_\alpha e_j := \alpha_j e_{j+1}$ $(j \geq 0)$ and extended by linearity. (It is standard that for our questions of interest we may, and do henceforth, assume that $\alpha$ is positive, i.e., $\alpha_n > 0$ for all $n$.) The **moments** $\gamma = (\gamma_n)_{n=0}^\infty$ of the shift are given by $\gamma_0 := 1$ and $\gamma_n := \prod_{j=0}^{n-1} \alpha_j^2$ for $n \geq 1$.

An operator $T \in \mathcal{L}(\mathcal{H})$ is normal if $TT^* = T^*T$, where $T^*$ is the adjoint of $T$. (It is well known that with our assumption of positive weights no weighted shift can be normal.) An operator is subnormal if it is the restriction of a normal operator to a (closed) invariant subspace. Weighted shifts have been studied extensively with respect to subnormality and related properties; there is a condition sufficient for subnormality via an approach through $k$-hyponormality (positive (semi)-definiteness of moment matrices) and the Bram-Halmos theorem (see originally [9], and [5] for background and further references). There is a second approach through $n$-contractivity and the Agler-Embry theorem (see originally [1] and again [5] for context); as well, a weighted shift is subnormal if and only if it has a Berger measure (a measure compactly supported in $\mathbb{R}_+$ whose moment sequence is the moment sequence of the shift – see [10, III.8.16] and [11]).

It is well known that the Schur (Hadamard, or entry-wise) product of two positive (semi)-definite matrices is again positive (semi)-definite. It follows readily from the Bram-Halmos approach to subnormality that if $W_\alpha$ and $W_\beta$ are subnormal (respectively, $k$-hyponormal for some $k$) shifts, then the shift with weight sequence the Schur product $\alpha\beta$ is subnormal (respectively, $k$-hyponormal).

In [3] and [4] certain "better than subnormal" shifts were considered. We say that a shift $W_\alpha$ is **moment infinitely divisible** (henceforth, $\mathcal{MID}$) if, for every real $s$ such that $s > 0$, the shift with weight sequence $\alpha^{(s)} = (\alpha_n^s)_{n=0}^\infty$ is subnormal. Recall that the forward difference operator $\nabla$ on a sequence $\beta = (\beta_n)$ is defined by $(\nabla\beta)_n := \beta_n - \beta_{n+1}$, and powers of $\nabla$ are defined recursively by $\nabla^{(0)}$ as the identity mapping and $\nabla^{(n+1)} := \nabla(\nabla^{(n)})$. A sequence $\beta$ is completely monotone if $\nabla^{(n)}\beta \geq 0$ for all $n = 1, 2, \ldots$. We say a sequence is **log completely monotone** if $\nabla^{(n)} \ln \beta \geq 0$ for all $n = 1, 2, \ldots$; note that we do not require that the terms $\ln \beta_n$ be themselves positive.

The backward difference operator $\Delta$ is defined by $\Delta := -\nabla$, and a sequence $\beta$ is **completely alternating** (respectively, **log completely alternating**) if $\Delta^{(n)}\beta \geq 0$ for $n = 1, 2, \ldots$ (respectively, $\Delta^{(n)}\ln\beta \geq 0$ for $n = 1, 2, \ldots$, where again, in contrast to the definition in [8] we do not require that the terms $(\ln\beta)_n$ be themselves positive).

It is known from [3] and [4] that a contractive shift is $\mathcal{MID}$ if and only if its moment sequence is log completely monotone, or, equivalently, its weight sequence is log completely alternating. These conditions may sometimes be verified by showing that the sequence is interpolated by a function in a certain class (Bernstein or log Bernstein); see [5] for definitions and discussion. Finally, an operator class related to subnormality is the class of completely hyperexpansive operators; these are obtained by reversing the direction of the inequality in the tests for $n$-contractivity (see originally [2], and [5] for definitions and discussion).

Recall from [5] that the geometrically regular weighted shifts (briefly, GRWS) are those shifts with weight sequence given by

$$\alpha_n(N, D) = \sqrt{\frac{p^n + N}{p^n + D}},$$

where $p > 1$ and $(N, D)$ is fixed in $(-1, 1) \times (-1, 1)$.

We state the needed partial versions of previous results for geometrically regular weighted shifts (themselves, as opposed to their quotients which are the subject of the current paper) from [5], and in particular their membership in the various general classes, in reference to the diagram in Figure 1. This is, for fixed $p > 1$, the open unit square in $\mathbb{R}^2$ with parameter pair $(N, D)$, and with Sectors I, II, $\ldots$, the ray $D = pN$ in Sector VIII and the rays of the form $D = p^n N$ in Sector IV. (Hereafter, by a sector we mean the convex cone in the open unit square with vertex at $(0, 0)$ and bounded by two rays emanating from the center. All sectors in this paper will be closed in the relative topology of the open unit square.) We take the slight liberty of referring to the open unit square in this context as the "magic square."

The classes of interest relevant here appear as follows:

- In Sector I, the weighted shifts are $\mathcal{MID}$ with the sequence of weights squared interpolated by a Bernstein function;
- In Sector II, the weighted shifts are $\mathcal{MID}$ with the sequence of weights squared interpolated by a log Bernstein function;
- In Sector III, the weighted shifts are subnormal;
- In Sector IV, the shifts corresponding to points on the special lines $D = p^n N$, $n = 0, 1, 2, \ldots$, are subnormal with finitely atomic Berger measure but not $\mathcal{MID}$ (except on the line $D = N$, which yields the unweighted shift, which is $\mathcal{MID}$); on the boundary with Sector III the shift is subnormal but not $\mathcal{MID}$;

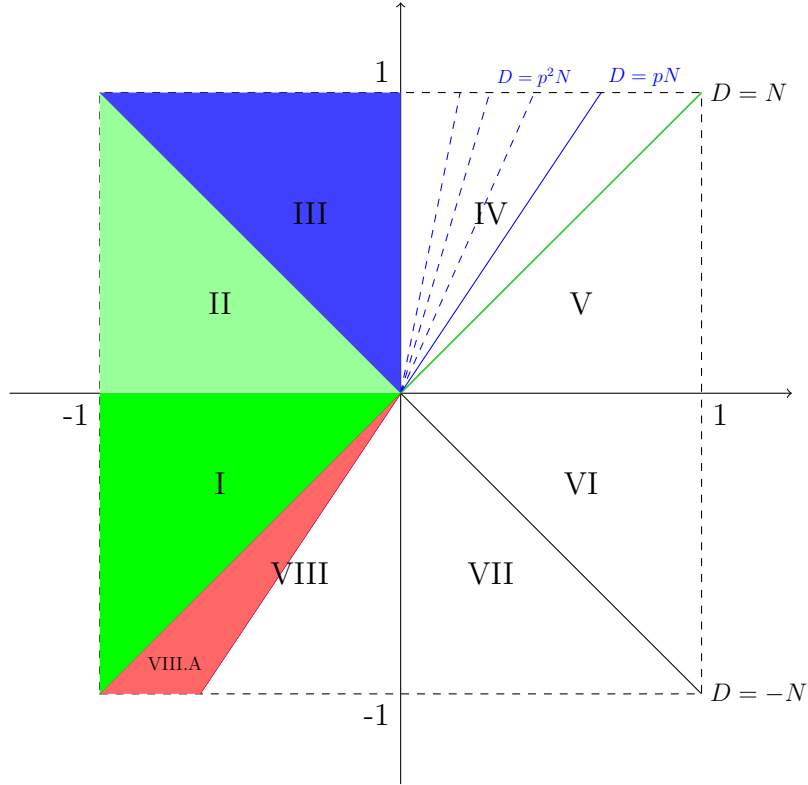

Figure 1. Magic Square

other points in the sector yield shifts that are not subnormal (although they are $k$-hyponormal for various $k$);
• In Sector VIIIA the shifts are completely hyperexpansive.

This numbering is chosen so that Sector I holds the "best" parameters, yielding shifts that are not only subnormal but $\mathcal{MID}$; indeed, their weights have the stronger property that they are interpolated by a Bernstein function, and, for a portion of the sector, the even stronger property of being the reciprocals of the weights of a completely hyperexpansive weighted shift (see [2] and [3] for the discussion showing this is stronger than $\mathcal{MID}$). Sector II contains the next best parameters, yielding an $\mathcal{MID}$ shift, and from there we simply proceed clockwise.

## 2. Main results

We first give a useful definition of two partial orders on the parameter space yielding the shifts we study.

**Definition 2.1.** Let $(N, D)$ and $(M, P)$ be two points in the open unit square $(-1, 1) \times (-1, 1)$. We say $(N, D) \gg (M, P)$ if the weight $\frac{\alpha(N,D)}{\alpha(M,P)}$ is an $\mathcal{MID}$ weight (i.e., the *shift with weights the Schur quotient is $\mathcal{MID}$). Often, particularly in text, we will say $(M, P)$ is $\mathcal{MID}$-subordinate to $(N, D)$ if $(N, D) \gg (M, P)$.*

Similarly, we say

**Definition 2.2.** $(N, D) \gg_s (M, P)$ *if the weight $\frac{\alpha(N,D)}{\alpha(M,P)}$ is a subnormal weight (i.e., the shift with weights the Schur quotient is subnormal). We will say $(M, P)$ is subnormal-subordinate to $(N, D)$ if $(N, D) \gg_s (M, P)$.*

Implicit in the above definitions, but omitted from the notation, is that we consider $p$ as fixed in such situations.

**Proposition 2.3.** *The relation "$\gg$" [respectively, "$\gg_s$"] is a partial order.*

Proof. The proof is straightforward using, for antisymmetry, that a subnormal shift must have a weakly increasing weight sequence, and, to handle transitivity, that the Schur product of two $\mathcal{MID}$ (respectively, subnormal) shifts is $\mathcal{MID}$ (respectively, subnormal). □

For a given "base point" $(N, D)$, define as well the sets

$$(2.1) \qquad \mathcal{MQ}_{(N,D)} := \{(M, P) \in (-1, 1) \times (-1, 1) \text{ such that } (N, D) \gg (M, P)\},$$

and

$$(2.2) \qquad \mathcal{SQ}_{(N,D)} := \{(M, P) \in (-1, 1) \times (-1, 1) \text{ such that } (N, D) \gg_s (M, P)\}.$$

We may refer to these sets as containing the points yielding "safe $\mathcal{MID}$ [respectively, subnormal] quotients."

Observe also that $\mathcal{MQ}_{(N,D)} \subseteq \mathcal{SQ}_{(N,D)}$.

Given a base point $(N, D)$ in one of various sectors, we will give sets contained in the sets $\mathcal{MQ}_{(N,D)}$ and $\mathcal{SQ}_{(N,D)}$; we will often refer to these as safe zones. The pictures in Figures 2, 3 and 4 contain the most important relevant results.

**Theorem 2.4.** *For $(N, D)$ in Sector I, $\mathcal{MQ}_{(N,D)}$ – the set of safe $\mathcal{MID}$ quotients – contains the set shown in Figure 2A; for $(N, D)$ in Sector II, $\mathcal{MQ}_{(N,D)}$ contains the set shown in Figure 2B.*

Similar results for base points in Sectors III-VIII and their $\mathcal{MID}$ safe quotients are given in the next section, where the results are proved.

**Theorem 2.5.** *For $(N, D)$ in Sector I, $\mathcal{SQ}_{(N,D)}$ – the set of safe subnormal quotients – contains the set shown in Figure 3A; for $(N, D)$ in Sector II, $\mathcal{SQ}_{(N,D)}$ contains the set shown in Figure 3B; for $(N, D)$ in Sector III, $\mathcal{SQ}_{(N,D)}$ contains the set shown in Figure 4.*

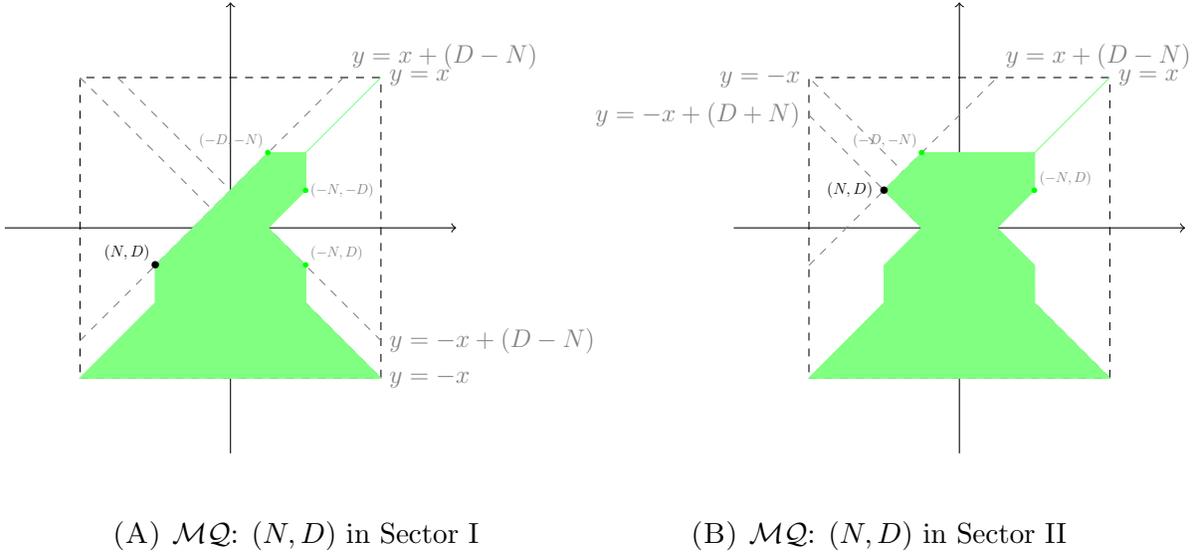

(A) $\mathcal{MQ}$: $(N, D)$ in Sector I

(B) $\mathcal{MQ}$: $(N, D)$ in Sector II

FIGURE 2. $\mathcal{MQ}$: $(N, D)$ in Sectors I or II

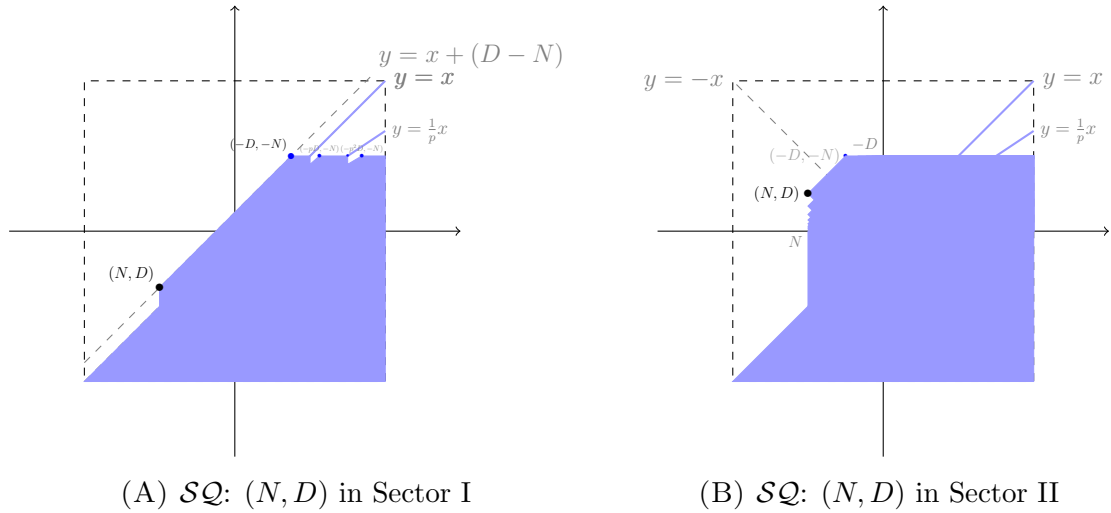

(A) $\mathcal{SQ}$: $(N, D)$ in Sector I

(B) $\mathcal{SQ}$: $(N, D)$ in Sector II

FIGURE 3. $\mathcal{SQ}$: $(N, D)$ in Sectors I or II

Note that the serrated left boundaries for $(N, D)$ in Sectors II and III are intentional.

Similar results for base points in the Sectors IV-VIII and their subnormal safe quotients are also given in the next section.

**Remark 2.6.** "Reflection principle"

We pause here to record an observation useful in the sequel. Observe that the weight sequence associated with a point $(M, P)$ in the magic square is the reciprocal of the weight sequence associated with the point $(P, M)$ (the reflection across the main diagonal $y = x$). In particular, and since we will be using Schur products throughout, the statement that "the quotient by $(M, P)$ is good for some property" is equivalent to "the product by $(P, M)$ is good for the same property." For example, in considering some quotient like $\alpha(N, D)/\alpha(M, P)$ for the $\mathcal{MID}$ property, we are well on our way to an $\mathcal{MID}$ shift if $\alpha(N, D)$ is $\mathcal{MID}$, since $\mathcal{MID}$ is preserved under Schur products. This will often provide zones of safe quotients "for free."

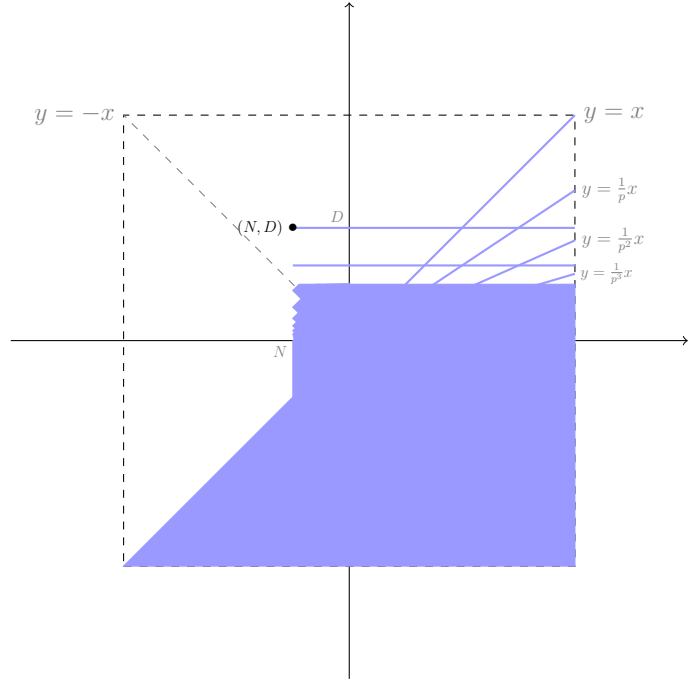

FIGURE 4. $\mathcal{SQ}$: $(N, D)$ in Sector III

This same observation allows for an immediate translation from the quotient results in the figures above to product results, and we will leave such translation statements to the interested reader.

## 3. Proofs of Results

The basic method of proof, which we will use throughout, is to express some quotient of weights as a Schur product of weights known to be in some class of interest. We will also use strongly the transitivity of the relations "$\gg$" and "$\gg_s$". Let us separate out four arguments which will be recurring.

### 3.1. Recurring Arguments.

(i) If $(N, D)$ is our base point, analysis of $(M, D)$ on either a vertical or horizontal line including $(N, D)$ is in general straightforward, because cancellation results in the weights for a GRWS, which we understand reasonably well, per the diagram in Figure 1. This argument is to be found, for example, in consideration of the "vertical segment down" from $(N, D)$ in Sector I, as in Section 4.2.

*(ii)* Particularly for diagonal segments, using an argument involving completely monotone families from [5], we will express some quotient of weights as interpolated by a Bernstein function, thus establishing that the weights correspond to an $\mathcal{MID}$ shift. This argument can be seen, for example, in the "diagonal segments over" $(N, D)$ in Sector I, as in Section 4.3.

*(iii)* We will use a simple but important application of transitivity. We may for some $(N, D)$ find a particularly useful point $(N', D')$ so that $(N', D')$ is $\mathcal{MID}$-subordinate to $(N, D)$ (respectively, $(N', D')$ is subnormal-subordinate to $(N', D')$) thus acquiring for $(N, D)$ the safe zone pertaining to $(N', D')$. This allows us to leverage prior results; an example here is the observation, for a base point $(N, D)$ in Sector IV, that we have $(N, D) \gg (-N, -D)$, as in Section 8.1.

*(iv)* Finally, we note that trivially either sort of safe quotient zone includes the base point $(N, D)$ itself. For clarity in diagrams, we will indicate $(N, D)$ by a black dot without giving it the color to which it is entitled.

The following statements are very useful in the sequel, as they provide methods for "growing" the zone of safe quotients. In what follows, some designations such as "NW-SE" refer to the compass directions "North West" and "South East."

**Theorem 3.1.** *For a given $(N, D) \in (-1, 1) \times (-1, 1)$, we have*

(1) $(M, P) \in \mathcal{MQ}_{(N,D)}$ *and* $-1 < M \leq 0$ *implies* $(-M, P) \in \mathcal{MQ}_{(N,D)}$. *(Right shadow, shown in )*

(2) $(M, P) \in \mathcal{MQ}_{(N,D)}$ *and* $|M| \geq |P|$ *implies* $(-P, -M) \in \mathcal{MQ}_{(N,D)}$. *(Central reversed shadow, shown in )*

(3)  a) $(M, P) \in \mathcal{MQ}_{(N,D)} \cap (-1, 0] \times [0, 1)$ *implies* $(-M, -P) \in \mathcal{MQ}_{(N,D)}$. *(NW-SE shadow, shown in )*

  b) $(M, P) \in \mathcal{MQ}_{(N,D)} \cap (-1, 0] \times (-1, 0]$ *and* $P \geq M$ *together imply* $(-M, -P) \in \mathcal{MQ}_{(N,D)}$. *(SW-NE shadow, shown in )*

Proof. 1) Assume $(M, P) \in \mathcal{MQ}_{(N,D)}$ and consider

$$\frac{\alpha(N, D)}{\alpha(-M, P)} = \frac{\alpha(N, D)}{\alpha(M, P)} \cdot \alpha(M, -M).$$

The weight $\alpha(M, -M)$ is obviously a $\mathcal{MID}$ weight since it corresponds to a GRWS with the pair living on the line $y = -x$ and $x \leq 0$ (cf. Figure 1). Since $\mathcal{MID}$ is stable under Schur products, we have as desired $(N, D) \gg (-M, P)$, i.e., $(-M, P) \in \mathcal{MQ}_{(N,D)}$ .

2) In the same manner, we write

$$\frac{\alpha(N,D)}{\alpha(-P,-M)} = \frac{\alpha(N,D)}{\alpha(M,P)} \cdot \frac{\alpha(M,P)}{\alpha(-P,-M)}.$$

Surely $(N,D) \gg (M,P)$ by definition if we assume $(M,P) \in \mathcal{MQ}_{(N,D)}$. With the assumption $|M| \geq |P|$, we have also $(-P,-M)$ is $\mathcal{MID}$-subordinate to $(M,P)$. Indeed,

$$\frac{\alpha(M,P)}{\alpha(-P,-M)} = \frac{p^{2n}-M^2}{p^{2n}-P^2} = \frac{(p^2)^n - M^2}{(p^2)^n - P^2}$$

and this corresponds to the GRWS associated to $p^2$ with the pair $(-M^2, -P^2)$ located in the associated Sector I. Thus by transitivity $(N,D) \gg (M,P) \gg (-P,-M)$ and we have the desired result.

3) a) As above

$$\frac{\alpha(N,D)}{\alpha(-M,-P)} = \frac{\alpha(N,D)}{\alpha(M,P)} \cdot \frac{\alpha(M,P)}{\alpha(-M,-P)} = \frac{\alpha(N,D)}{\alpha(M,P)} \cdot \alpha(M,-M) \cdot \alpha(-P,P)$$

and the result follows since we have the Schur product of three $\mathcal{MID}$ weights.

For 3) b) we need the following lemma.

**Lemma 3.2.** *Let $(M,P)$ be a pair in $(-1,1) \times (-1,1)$ such that $P \geq M$ and $MP \geq 0$; then*

$$(M,P) \gg (-M,-P).$$

**Proof.** Notice that

$$\frac{\alpha(M,P)}{\alpha(-M,-P)} = \frac{(p^n+M)(p^n-P)}{(p^n-M)(p^n+P)},$$

so if $P = M$ the result is trivial. Otherwise, we would like to prove that the interpolating function, namely

$$g(x) := \frac{(p^x+M)(p^x-P)}{(p^x-M)(p^x+P)},$$

is a log Bernstein function. (The reader of [3] may expect that this function should be contractive; while we note that it is, it is now understood that this restriction is not necessary.)

First, it is easy to see that $0 < g(x) < 1$ for all $x > 0$. Second, we have

$$(\log g)'(x) = \frac{g'(x)}{g(x)} = \frac{2(\log p)(P-M)p^x(p^{2x}+MP)}{(p^{2x}-M^2)(p^{2x}-P^2)}.$$

Observe that this is positive as required. Clearly we may discard the $2(P-M)$ and the $\log p$. We then are left with

$$\frac{p^x(p^{2x}+MP)}{(p^{2x}-M^2)(p^{2x}-P^2)} = \frac{p^x}{p^{2x}-P^2} \cdot \frac{p^{2x}}{p^{2x}-M^2} + \frac{MPp^x}{p^{2x}-M^2} \cdot \frac{1}{p^{2x}-P^2}$$

which is completely monotone using the result on completely monotone families ([5, Lemma 2.2]) with the specific family

$$\left\{ \frac{1}{p^{2x} - M^2} \ , \ \frac{p^x}{p^{2x} - M^2} \ , \ \frac{p^{2x}}{p^{2x} - M^2} \ , \ \frac{1}{p^{2x} - P^2} \ , \ \frac{p^x}{p^{2x} - P^2} \ , \ \frac{p^{2x}}{p^{2x} - P^2} \right\}.$$

Now, to prove (3) b), let $(M, P)$ be in $\mathcal{MQ}_{(N,D)} \cap (-1, 0] \times (-1, 0]$ and $P \geq M$ and just use the above lemma and transitivity: $(N, D) \gg (M, P) \gg (-M, -P)$. $\qquad \square$

With these results in hand, we turn to producing the various safe zones.

## 4. Safe $\mathcal{MID}$ quotients for a base point in Sector I

Let $(N, D)$ be a pair in Sector I; that is, $0 \geq D \geq N > -1$. By the "Reflection Principle" it follows readily that for a point $(N, D)$ in either Sector I or Sector II, Sectors VII and VIII provide safe quotients, and so these latter sectors are automatically included in the zone for safe quotients.

Note second that the main diagonal in the magic square yields simply the unilateral shift (which is both $\mathcal{MID}$ and completely hyperexpansive) and thus any point on this diagonal yields a safe quotient for any property of interest which the base point $(N, D)$ itself possesses.

We may now turn to further zones.

4.1. **The green zone.** The aim is to prove that the green zone in Figure 5 consists of points $(M, P)$ which are $\mathcal{MID}$-subordinate to $(N, D)$ and thus yielding safe $\mathcal{MID}$ quotients.

We will proceed in three steps.

4.2. **Vertical segment down from** $(N, D)$**.** Suppose that $(M, P)$ is on the green segment as in Figure 6A, which means that $(M, P) = (N, P)$ where $N \leq P \leq D$. It is clear that

$$\frac{\alpha(N, D)}{\alpha(M, P)} = \frac{\alpha(N, D)}{\alpha(N, P)} = \alpha(P, D)$$

which is obviously an $\mathcal{MID}$ weight since the pair $(P, D)$ is in Sector I.

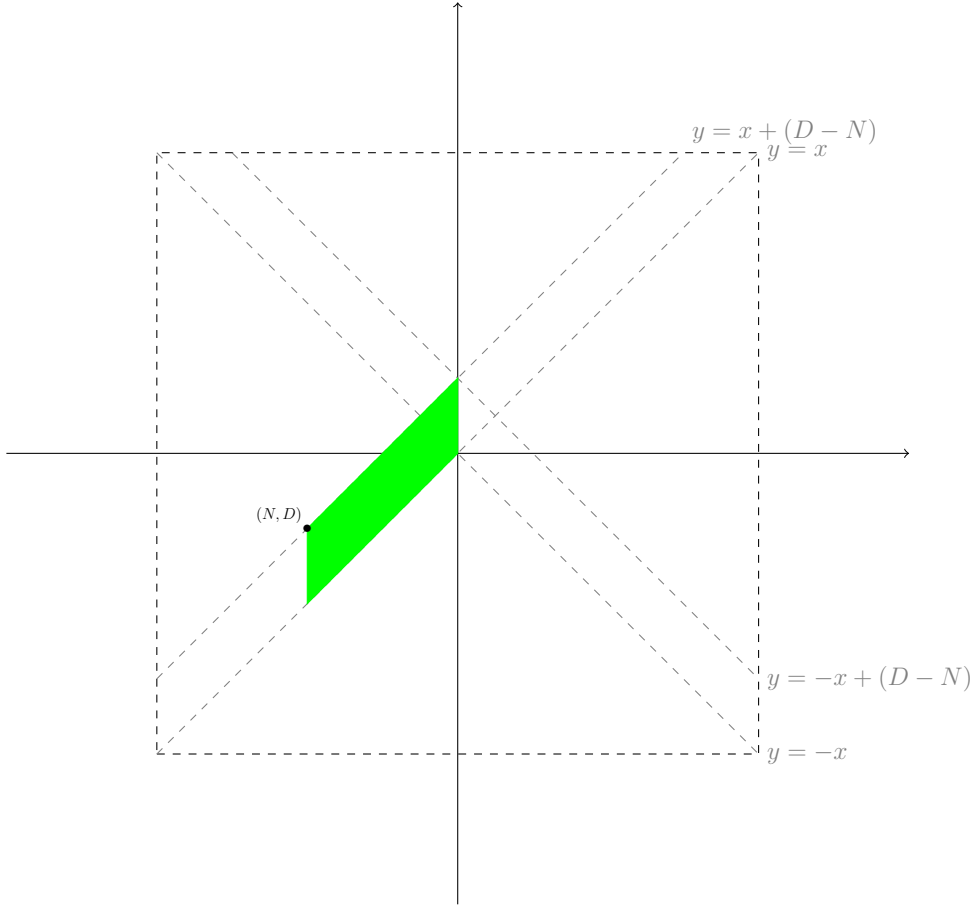

FIGURE 5. The green zone

4.3. **Diagonal segment over** $(N, D)$**.** Suppose that $(M, P)$ is in the diagonal green segment as in Figure 6A, so $(M, P) = (N, D) + q(1, 1)$ and $0 \leq q \leq -N$. Then

$$\begin{aligned}
\frac{\alpha^2(N, D)}{\alpha^2(M, P)} &= \frac{\alpha^2(N, D)}{\alpha^2(N + q, D + q)} \\
&= \frac{p^n + N}{p^n + D} \cdot \frac{p^n + D + q}{p^n + N + q} \\
&= 1 - \frac{q(D - N)}{(p^n + D)(p^n + N + q)}
\end{aligned}$$

We claim that $x \mapsto \dfrac{q(D - N)}{(p^x + D)(p^x + N + q)}$ is a contractive completely monotone function; this uses the result on completely monotone families ([5, Lemma 2.2]) where the involved family is

$$\left\{ \frac{1}{p^x + D}, \frac{p^x}{p^x + D}, \frac{1}{p^x + N + q}, \frac{p^x}{p^x + N + q} \right\}, \text{ with } -1 < D < 0 \text{ and } -1 < N + q < 0.$$

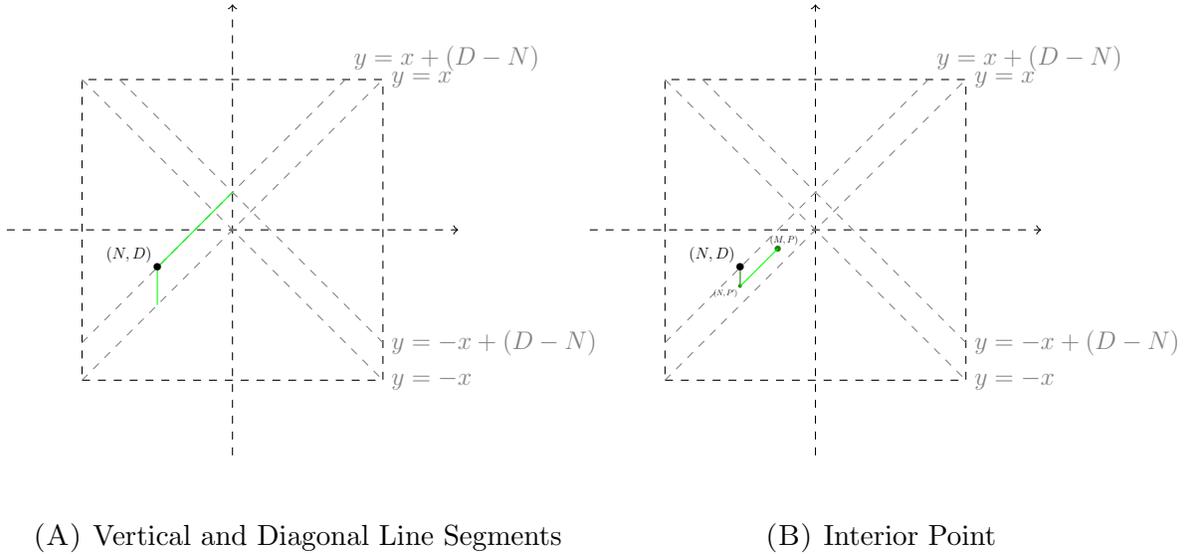

(A) Vertical and Diagonal Line Segments          (B) Interior Point

FIGURE 6. Work toward green zone

Therefore, the function $f(x) := 1 - \dfrac{q(D-N)}{(p^x + D)(p^x + N + q)}$ is a Bernstein function. This insures that $\dfrac{\alpha(N, D)}{\alpha(N+q, D+q)}$ is a $\mathcal{MID}$ weight.

Therefore,

$$(N, D) \gg (N+q, D+q) \quad \text{for every} \quad 0 \le q \le -N.$$

**Remark 4.1.** Notice here that according to step 2, for any given point $(N', D')$ in Sector I, the analogous diagonal segment up and to the right from it consists of points $\mathcal{MID}$-subordinate to $(N', D')$.

4.4. **The last step.** For the general situation (when the pair $(M, P)$ is in the interior of the green parallelogram), we may use the path $(N, D) \longrightarrow (N, P') \longrightarrow (M, P)$ as in Figure 6B, to prove our result. Indeed, $(N, D) \gg (N, P') \gg (M, P)$, and this yields the result. (As we have previously done in our diagrams, we will continue to indicate $(N, D)$ by a black dot without giving it the color to which it is entitled.)          $\square$

4.5. **Central reversed shadow (The blue zone).** Using the second statement of Theorem 3.1, we obtain the blue zone in Figure 7A (notice that this is nothing but symmetry with respect to the axis $y = -x$):

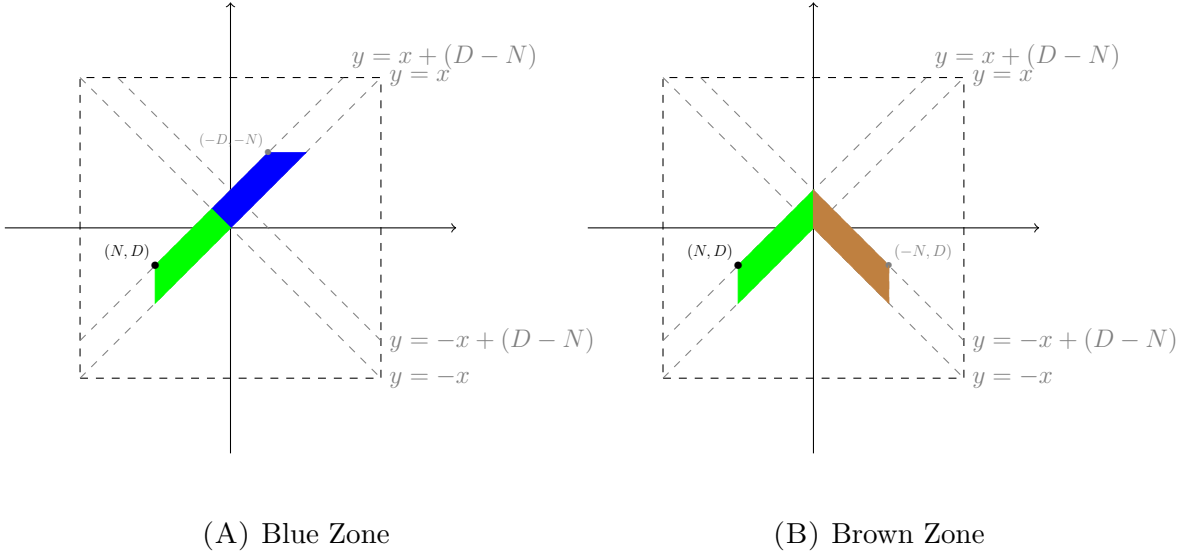

(A) Blue Zone                    (B) Brown Zone

Figure 7. Blue and Brown Zones

4.6. **Right shadow (The brown zone).** Using Theorem 3.1 (1), we obtain the brown zone (cf. Figure 7B).

**Remark 4.2.** Note that the purple zone (NW-SE shadow) – a small triangle in the fourth quadrant arising from a small triangle in the second quadrant – is contained in, and pictorially hidden by, the brown zone.

4.7. **SW-NE shadow (The red zone).** Using Theorem 3.1 (3) b), we obtain this zone, as in Figure 8A.

4.8. **Conclusion.** By gathering the above zones and including Sectors VII and VIII as noted at the beginning of the discussion for a base point in Sector I, we obtain the safe zone for $\mathcal{MID}$ quotients as in Figure 8B or as in Figure 2A in Section 2 (Main Results) on page 6.

## 5. Safe $\mathcal{MID}$ quotients for a base point in Sector II

Our assumption is $-N \geq D > 0 > N > -1$. We begin with the following result.

**Theorem 5.1.** *Let $(N, D)$ be in Sector II ($-N \geq D > 0 > N > -1$). Then we have*

$$\mathcal{MQ}_{(N,-D)} \subseteq \mathcal{MQ}_{(N,D)}.$$

Proof.   Just write

$$\frac{\alpha(N, D)}{\alpha(M, P)} = \frac{\alpha(N, D)}{\alpha(N, -D)} \cdot \frac{\alpha(N, -D)}{\alpha(M, P)} = \alpha(-D, D) \cdot \frac{\alpha(N, -D)}{\alpha(M, P)}.$$

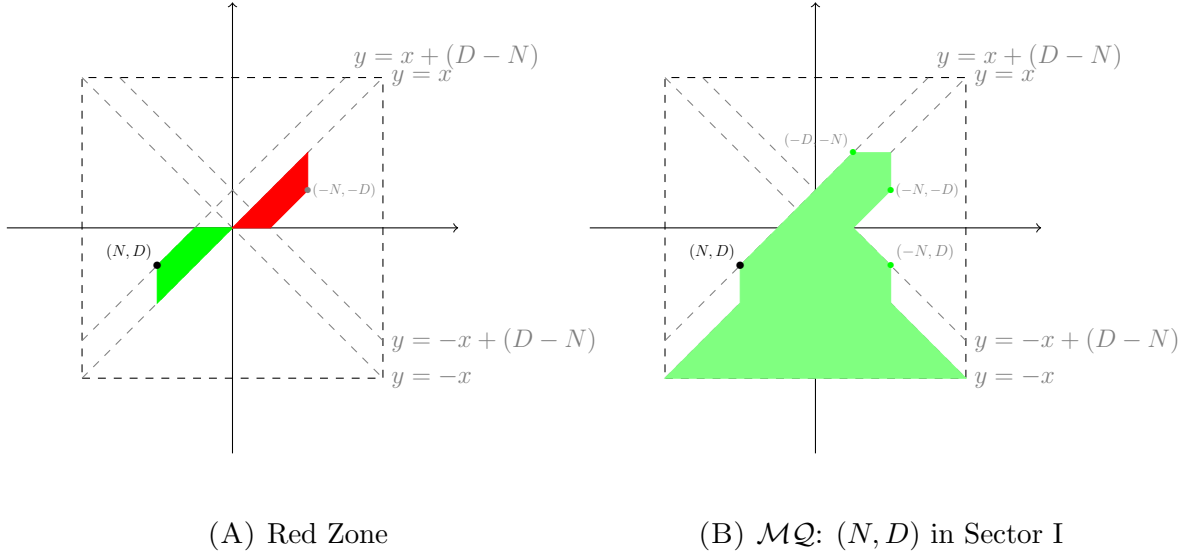

(A) Red Zone                    (B) $\mathcal{MQ}$: $(N, D)$ in Sector I

FIGURE 8. Red Zone and $\mathcal{MQ}$: $(N, D)$ in Sector I

Since $(-D, D)$ is in Sector II, $\alpha(-D, D)$ is an $\mathcal{MID}$ weight, and we have the result. $\square$

**Remark 5.2.** The application of this theorem is straightforward and we will reserve its implications until we assemble results into a picture at the end. Observe that

$$(N, D) \gg (M, P) \iff (N, M) \gg (D, P).$$

The most useful approach for this sector is to take advantage of the known facts about safe zones for the quotients with a starting point in Sector I.

5.1. **"Left side" (The green zone).** If $M < 0$ then $(N, M)$ is in Sector I if furthermore $N \leq M$. Thus, in this case (since we know already that $-N \geq D > 0$), we have

$$(N, M) \gg (D, P) \iff \begin{cases} D - (M - N) \leq P \leq D + (M - N) \\ D \leq -N \text{ (which holds)} \\ 0 \leq P \leq -N \end{cases}$$
$$\iff \begin{cases} -M + (D + N) \leq P \leq M + (D - N) \\ 0 \leq P \leq -N \end{cases}$$

Thus we obtain the safe zone (Figure 9A, green zone).

5.2. **"Right shadow" (The brown zone).** Using the first statement of Theorem 3.1, we get (Figure 9A, brown zone):

5.3. **"Central reversed shadow" (The blue zone).** Using the second statement of Theorem 3.1, we obtain Figure 9B.

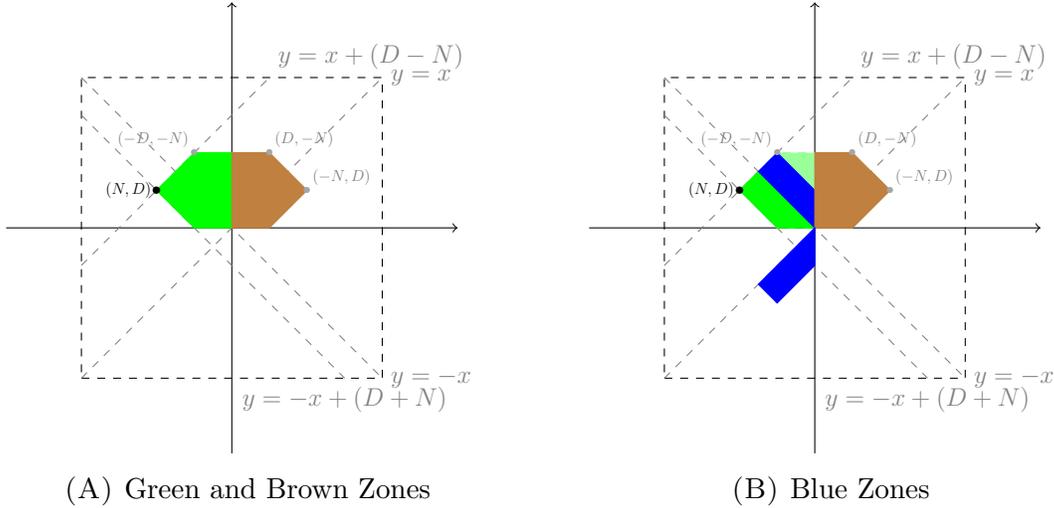

(A) Green and Brown Zones

(B) Blue Zones

FIGURE 9. Green, Brown, and Blue Zones

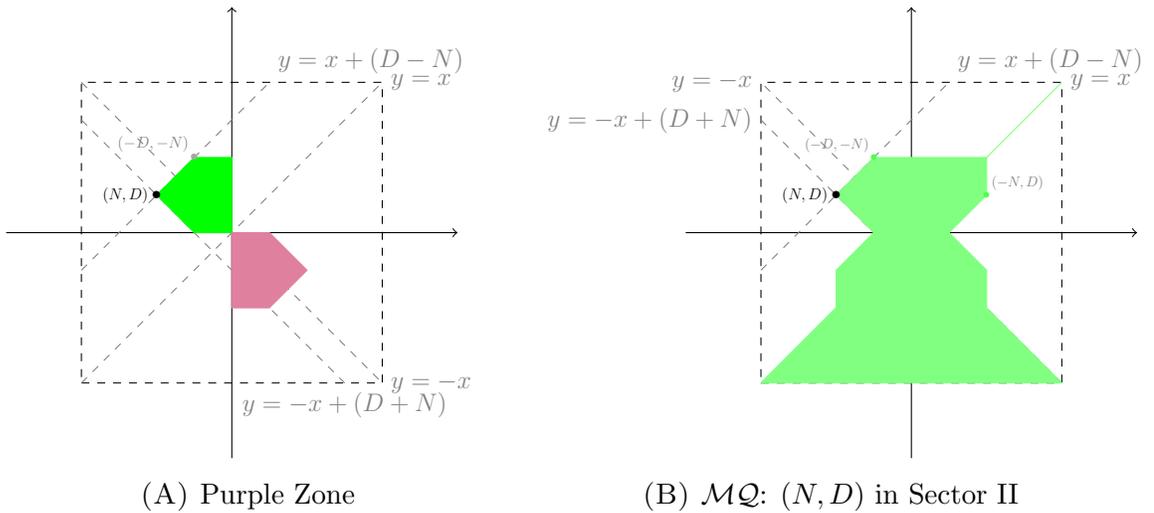

(A) Purple Zone

(B) $\mathcal{MQ}$: $(N, D)$ in Sector II

FIGURE 10. Purple Zone and $\mathcal{MQ}$: $(N, D)$ in Sector II

5.4. **"NW-SE shadow" (The purple zone).** Applying the third result of Theorem 3.1, we obtain the purple zone in Figure 10A.

**Remark 5.3.** Notice that the blue zone and the purple one are contained in, and in the pictures hidden by, other safe zones (see below).

5.5. **Conclusion.** Putting all of this together, including that Sectors VII and VIII are safe as noted before, and including the one resulting from Theorem 5.1, we obtain the picture in Figure 10B or as in Figure 2B in Section 2 (Main Results) on page 6.

We now turn to Sector VIII, departing from the order one might anticipate, so as to be able to leverage earlier results more efficiently.

## 6. Base point in Sector VIII

Consider $(N, D)$ in Sector VIII, which is $D < N < 0$, and some $(M, P)$ $\mathcal{MID}$-subordinate to $(N, D)$, so

$$\frac{\alpha(N, D)}{\alpha(M, P)} = \frac{p^n + N}{p^n + M} \cdot \frac{p^n + P}{p^n + D}$$

is an $\mathcal{MID}$ weight.

A sufficient condition, using the terms in the above choice of grouping terms, conditions for $\mathcal{MID}$ GRWS as in Figure 1 (page 6), and Schur products, is

$$\left\{ \begin{array}{l} N \leq M \leq -N \\ P \leq D \leq -P \end{array} \right. \iff \left\{ \begin{array}{l} N \leq M \leq -N \\ -1 < P \leq D \end{array} \right.$$

which gives the first part of the corresponding safe zone.    (NB: In particular, the vertical line segment down from $(N, D)$ is safe.)

Now, we'll use a technique we have used above for when the pair $(M, P)$ is on a diagonal from the base point, but this time below $(N, D)$. This means $(M, P) = (N, D) + q(1, 1)$ and $-1 - N < q < 0$. Then

$$\begin{array}{rl} \dfrac{\alpha^2(N, D)}{\alpha^2(M, P)} & = \dfrac{\alpha^2(N, D)}{\alpha^2(N + q, D + q)} \\[2mm] & = \dfrac{p^n + N}{p^n + D} \cdot \dfrac{p^n + D + q}{p^n + N + q} \\[2mm] & = 1 - \dfrac{q(D - N)}{(p^n + D)(p^n + N + q)} \end{array}$$

We have $x \mapsto \dfrac{q(D - N)}{(p^x + D)(p^x + N + q)}$ is a completely monotone function using completely monotone families from [5] with the family

$\left\{ \dfrac{1}{p^x + D}, \dfrac{p^x}{p^x + D}, \dfrac{1}{p^x + N + q}, \dfrac{p^x}{p^x + N + q} \right\}$ with $-1 < D < 0$ and $-1 < N + q < 0$.

Therefore the function $f$ defined by

$$f(x) := 1 - \frac{q(D - N)}{(p^x + D)(p^x + N + q)}$$

is a Bernstein function. This ensures that

$$\frac{\alpha(N, D)}{\alpha(N + q, D + q)}$$

is a $\mathcal{MID}$ weight. So $(N, D) \gg (N + q, D + q)$ for every $q$, $-1 - N < q < 0$.

Using transitivity, we obtain the picture of the $\mathcal{MID}$ safe zone as in Figure 11A.

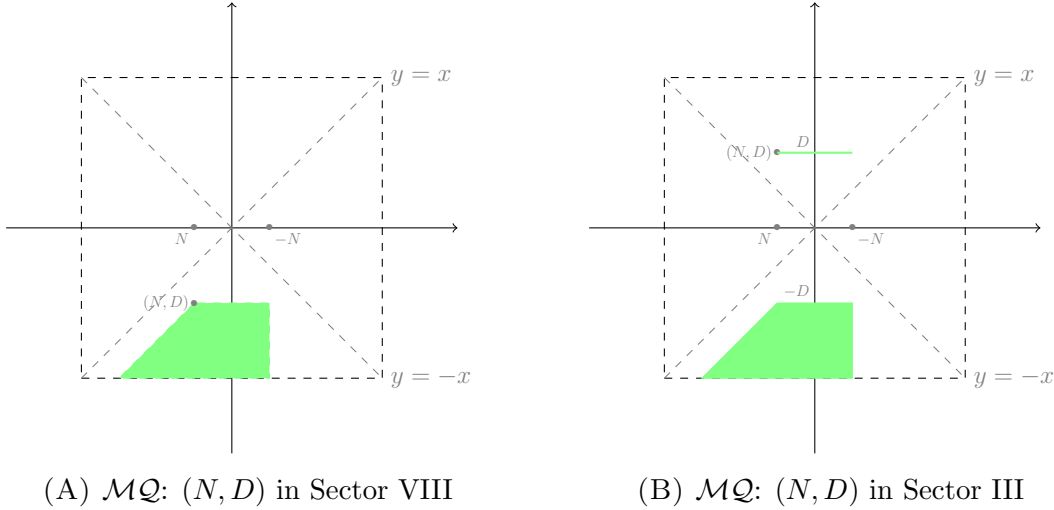

(A) $\mathcal{MQ}$: $(N, D)$ in Sector VIII

(B) $\mathcal{MQ}$: $(N, D)$ in Sector III

FIGURE 11. $\mathcal{MQ}$: $(N, D)$ in Sectors VIII or III

## 7. Base point in Sector III

Suppose $(N, D)$ is in Sector III, which is $N < 0$ and $D > -N$, and as usual consider some point $(M, P)$ subordinate to $(N, D)$. If $P = D$, it is clear that

$$\frac{\alpha(N, D)}{\alpha(M, P)} = \alpha(N, M)$$

is $\mathcal{MID}$ if $N \le M \le -N$ which proves that a certain horizontal line segment is safe.

Further, $\frac{\alpha(N,D)}{\alpha(N,-D)} = \alpha(-D, D)$, and this latter weight yields $\mathcal{MID}$ using $D > 0$ as it is in (on the boundary of) Sector II, so $\alpha(N, D) \gg \alpha(N, -D)$. With $(N, D)$ in Sector III, $(N, -D)$ is in Sector VIII, and we may apply transitivity; combining with the result from the previous paragraph, we obtain the $\mathcal{MID}$ safe zone in Figure 11B.

## 8. Base points in Sector IV, V, VI and VII

First, if $(N, D)$ is in Sector V, so $0 < D < N$, we prove that $(-D, -N)$ is $\mathcal{MID}$-subordinate to $(N, D)$ and we obtain a safe zone using the result for Sector VIII above. Observe that

$$\frac{\alpha(N, D)}{\alpha(-D, -N)} = \frac{p^{2n} - N^2}{p^{2n} - D^2}$$

and this is obviously a weight corresponding to an $\mathcal{MID}$ GRWS with $p^2$ given our assumptions on $N$ and $D$, yielding the picture in Figure 12B.

Second, if $(N, D)$ is in sector IV, we have

$$(8.1) \qquad\qquad\qquad 0 < N < D,$$

and citing Lemma 3.2 we obtain $(-N, -D)$ is $\mathcal{MID}$-subordinate to $(N, D)$. Using again transitivity and prior results we obtain the picture in Figure 12A.

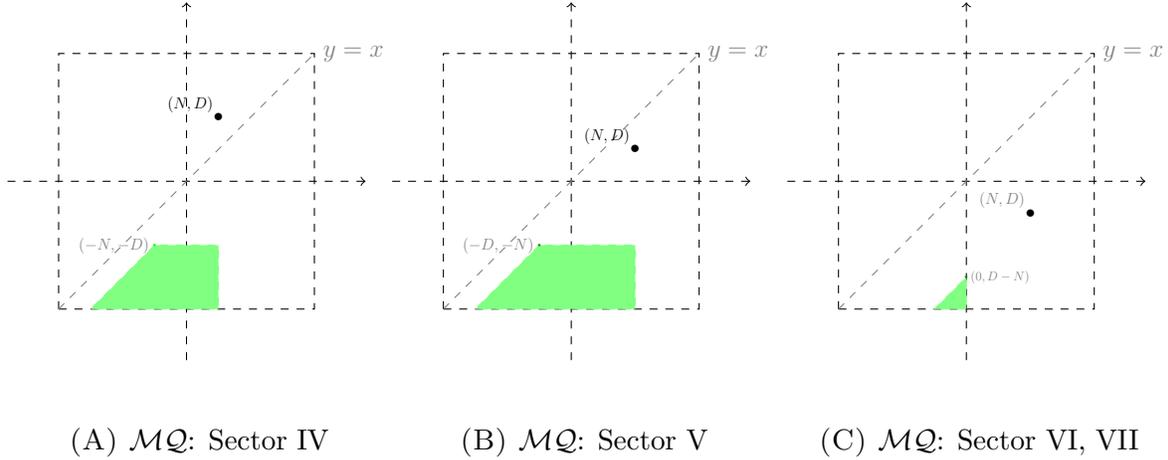

(A) $\mathcal{MQ}$: Sector IV        (B) $\mathcal{MQ}$: Sector V        (C) $\mathcal{MQ}$: Sector VI, VII

FIGURE 12. $\mathcal{MQ}$: $(N, D)$ in Sectors IV, V, VI, and VII

Third, if $(N, D)$ is in sector VI or VII, meaning that $D < 0 < N$, we'll show $(0, D - N)$ is $\mathcal{MID}$-subordinate to $(N, D)$. Indeed,

$$
\begin{aligned}
\frac{\alpha^2(N, D)}{\alpha^2(0, D - N)} &= \frac{p^n + N}{p^n + D} \cdot \frac{p^n + D - N}{p^n} \\
&= 1 - \frac{N(N - D)}{p^n(p^n + D)}.
\end{aligned}
$$

Since $-1 < D \le 0$ and $N(N - D) > 0$, the map $x \mapsto \dfrac{N(N - D)}{p^x(p^x + D)}$ is a completely monotone function using completely monotone families and so the function

$$
f(x) := 1 - \frac{N(N - D)}{p^x(p^x + D)}
$$

is a Bernstein function. Consequently,

$$
\frac{\alpha(N, D)}{\alpha(0, D - N)}
$$

is an $\mathcal{MID}$ weight so $(N, D) \gg (0, D - N)$, and again we apply transitivity. This gives the result in Figure 22 (which covers both Sectors VI and VII).

The three results obtained for $\mathcal{MID}$ safe quotients for a base point in Sector IV, V, and VI/VII are summarized in Figures 12A through 12C.

## 9. Subnormal safe quotients

We turn next to determining points yielding safe subnormal, as opposed to safe $\mathcal{MID}$, quotients. For $(N, D)$ and $(M, P)$ in our magic square $(-1, 1) \times (-1, 1)$, recall the partial order defined by:

$$
(N, D) \gg_s (M, P) \overset{\text{def}}{\Longleftrightarrow} \frac{\alpha(N, D)}{\alpha(M, P)} \quad \text{is a subnormal weight,}
$$

and, for a given pair $(N, D)$, the set

$$(9.1) \qquad \mathcal{SQ}_{(N,D)} := \{(M, P) \in (-1, 1) \times (-1, 1) \text{ such that } (N, D) \gg_s (M, P)\}.$$

**Remark 9.1.** 1) It is clear that $\mathcal{MQ}_{(N,D)} \subseteq \mathcal{SQ}_{(N,D)}$, which makes for a nesting of safe quotient zones.

2) The presence of the special lines in Sector IV, upon which there is subnormality, makes for considerably more complicated pictures for subnormal safe quotients than for $\mathcal{MID}$ quotients.

In the sequel, we will consider, for various placements of the base point, at least some portion of those zones. We again assemble some general observations.

First, the presence of the special lines creates what we call "raindrops" for base points $(N, D)$ with $D > 0$ (see Figure 14). If we consider points on the vertical line segment from $(N, D)$ with coordinates $(N, D/p^k)$ for some positive integer $k$, we have the quotient

$$\frac{\alpha(N, D)}{\alpha(N, D/p^k)_n} = \frac{\frac{p^n + N}{p^n + D}}{\frac{p^n + N}{p^n + D/p^k}} = \frac{p^n + D/p^k}{p^n + D}.$$

But this is a GRWS in Sector IV on the special line $y = p^k x$, and thus yields a subnormal shift. A similar result concerning "raindrops," but this time lying horizontally from $(N, D)$ and at $(p^k N, D)$, obtains if $N$ is positive.

Second, note also that a point $(N, P)$ on this same vertical line, but with $P \leq 0$, yields, by a similar computation for the quotient, a point in Sector II or III, again yielding a subnormal shift. A similar result for a horizontal line holds in the case in which $N$ is negative. These two yield, in what follows, horizontal or vertical line segments that extend "to the boundary" of the magic square.

Third, if $(N, D)$ produces raindrops of either kind, by transitivity the safe zone for subnormal quotients of $(N, D)$ will include the union of all the safe zones of the raindrops; it is this that produces the serrated boundaries for base points $(N, D)$ in Sectors II and III.

Fourth, if the base point $(N, D)$ already yields a subnormal (or $\mathcal{MID}$) shift, points on the reflections into Sector V of the special lines in Sector IV produce, via the reflection principle, subnormal shifts when used to form quotients.

### 9.1. $(N, D)$ **in Sector VIII.**

In this case, of course we contain the safe $\mathcal{MID}$ quotient zone as in Figure 11A; the picture in Figure 13A adds a zone arising from "horizontal lines that extend to the boundary" from points in the $\mathcal{MID}$ quotient zone as in the third observation above, and we leave the computation to the interested reader.

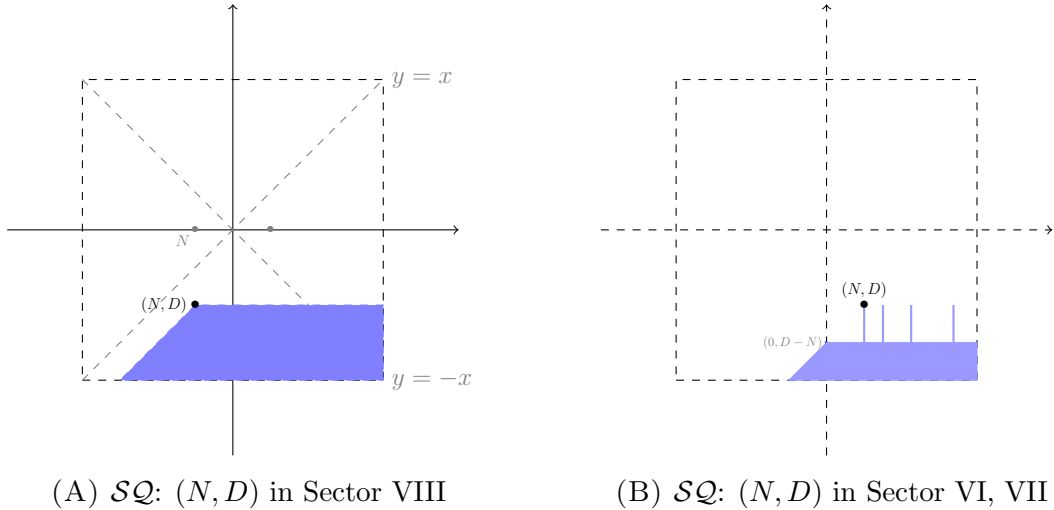

(A) $\mathcal{SQ}$: $(N, D)$ in Sector VIII

(B) $\mathcal{SQ}$: $(N, D)$ in Sector VI, VII

FIGURE 13. $\mathcal{SQ}$: $(N, D)$ in Sectors VIII, VI, or VII

**9.2. $(N, D)$ in Sectors VI or VII.** Here we have $N > 0$ and $D < 0$, and using the grouping

$$(9.2) \qquad \frac{\alpha(N, D)}{\alpha(M, P)} = \frac{p^n + N}{p^n + M} \cdot \frac{p^n + P}{p^n + D},$$

we have that if $M = p^k N$ for some integer $k$ and $-1 < P \leq D$ then $(N, D) \gg_s (M, P)$. Thus we obtain horizontal raindrops and vertical line segments that extend to the boundary from them. Further, one computes that $(0, D - N)$ is subnormal-subordinate to $(N, D)$; since this is a point in Sector VIII, we obtain by transitivity its safe subnormal quotient zone as just above, yielding the picture in Figure 13B.

**9.3. $(N, D)$ in Sector V.** We may obtain the horizontal and vertical raindrops as usual, as well as the aligned vertical lines down from the $x$ axis to the boundary. We have proved that for a point in Sector V we have $(-D, -N)$ is $\mathcal{MID}$-subordinate to $(N, D) \gg (-D, -N)$ (Section 8), so it is surely subnormal-subordinate, thus adding to the safe subnormal quotient zone for $(N, D)$ that for $(-D, -N)$ in Sector VIII, and Figure 14A follows.

**9.4. $(N, D)$ in Sector IV, "generic" case (not on a special line).** Here we have both $N > 0$ and $D > 0$, so we get both horizontal and vertical raindrops; as well, we obtain the vertical lines to the boundary in quadrant IV below the vertical raindrops. From the raindrops which lie in Sector V we obtain their safe zones as in Figure 14A from transitivity, and while no such raindrop can lie on the main diagonal $y = x$ (else $(N, D)$ would be on a special line) we obtain their union. In particular, we obtain the safe zone arising from the raindrop $(M, P) = (p^k N, D/p^j)$ in Sector V where $D/p^j$ is as large as possible and $p^k N$ is as small as possible; this is the first raindrop in the first

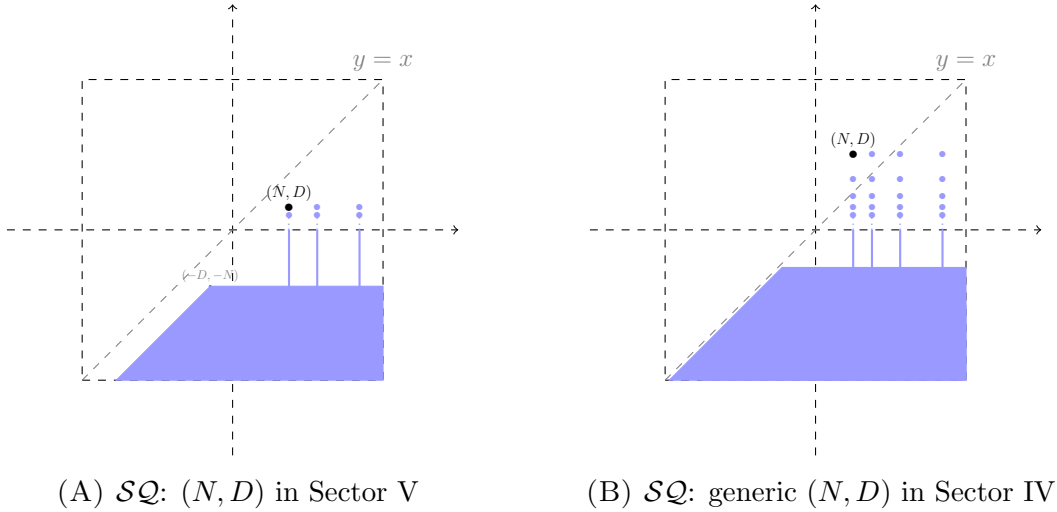

(A) $\mathcal{SQ}$: $(N, D)$ in Sector V

(B) $\mathcal{SQ}$: generic $(N, D)$ in Sector IV

FIGURE 14. $\mathcal{SQ}$: $(N, D)$ in Sector V or generic in Sector IV

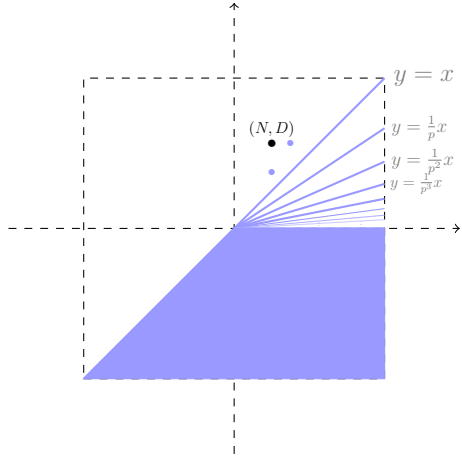

FIGURE 15. $\mathcal{SQ}$: non-generic $(N, D)$ in Sector IV

column which falls into Sector V. Putting this together yields the picture in Figure 14B as the safe subnormal quotient zone for a generic point in Sector IV.

9.5. $(N, D)$ **in Sector IV, "non-generic" case (on a special line).** In this case everything holding for the generic case still holds; as well, we obtain the reflections across the main diagonal into Sector V of the special lines in Sector IV (since the shift associated with $(N, D)$ is subnormal). Note that the raindrops from $(N, D)$ which lie in Sector V are actually on these reflected lines and are not otherwise indicated. We obtain all of Sectors VI, VII, and VIII by the reflection principle, since quotients with these points are products with points of subnormality. What results is the picture in Figure 15.

9.6. $(N, D)$ **in Sector I.** We begin with the safe zone for $\mathcal{MID}$ quotients as in Figure 2A on page 6. To this we may adjoin the reflections of special lines from Sector IV as usual; note that vertical or horizontal raindrops from points on these reflections will themselves be on reflected lines and add nothing new. Of course we obtain all of Sectors VI, VII, and VIII.

Consider next some horizontal segment in the safe zone and in the first quadrant. The Sector IV or Sector V first horizontal raindrops from the points in this horizontal segment create another horizontal segment shifted to the right; in particular, a segment from $(M, P)$ to $(P, P)$ is shifted to the segment from $(pM, P)$ to $(pP, P)$. These shifted segments will overlap if $pM \leq P$, thus extending the horizontal segment without a break. The raindrops from the second segment, of course, create a further shifted segment and so on, and it is easy to check that the overlaps continue if $pM \leq P$. What results is a horizontal segment, without breaks, from $(M, P)$ to the right-hand boundary.

We shall apply this general observation to the points on the left-hand upper diagonal boundary of the safe zone for $\mathcal{MID}$ quotients as in Figure 2A on page 6. The intersection of this boundary ($y = x + (D - N)$) with the line $y = px$ is at the point $(\frac{D-N}{p-1}, p\frac{D-N}{p-1})$. For any point $(Q, R)$ on the $y = x + (D - N)$ line with second coordinate less than $p\frac{D-N}{p-1}$), by the argument above we may add to the safe zone the whole horizontal line from $(Q, R)$ to the right-hand boundary of the square.

For points $(Q, R) = (Q, Q + (D - N))$ on the boundary line $y = x + (D - N)$ with $R = Q + (D - N) \geq -D$, the horizontal segment from $(Q, R)$ to $(-N, R)$ will overlap with its shifted segment $(pQ, R)$ to $(p(-N), R)$ if and only if $pQ \leq -N$, and this overlap of shifted segments will persist to create the unbroken horizontal line from $(Q, R)$ to the right-hand boundary of the square. If $pQ > -N$, however, the shifted segments will not overlap and we will create a "notch" with vertical right-hand boundary because of the common right shift of the vertical segment with first coordinate $-N$. The shifts of these segments may create a second notch, and so on.

There are a number of cases, and we content ourselves with drawing, in Figure 3A on page 6, a reasonably typical safe subnormal quotient zone for $(N, D)$ in Sector I.

9.7. $(N, D)$ **in Sector II.** We begin with the safe $\mathcal{MID}$ quotient zone as in Figure 2B on page 6. Since we have subnormality at $(N, D)$ we may add all of Sectors VI, VII, and VIII, and the reflections of the special lines from Sector IV into Sector V. Consideration (as in the discussion for Sector I) of horizontal segments from some point $(0, R)$ to $(R, R)$ with $0 \leq R \leq -N$ fill in the rectangle with corners $(0, 0)$, $(0, -N)$, $(1, -N)$, and $(1, 0)$. Consideration of the vertical segment from $(N, -D)$ to $(N, 0)$, easy as discussed before, yields that this segment is safe, and then a repetition of the diagonal argument from it yields (so far) the picture in Figure 16A.

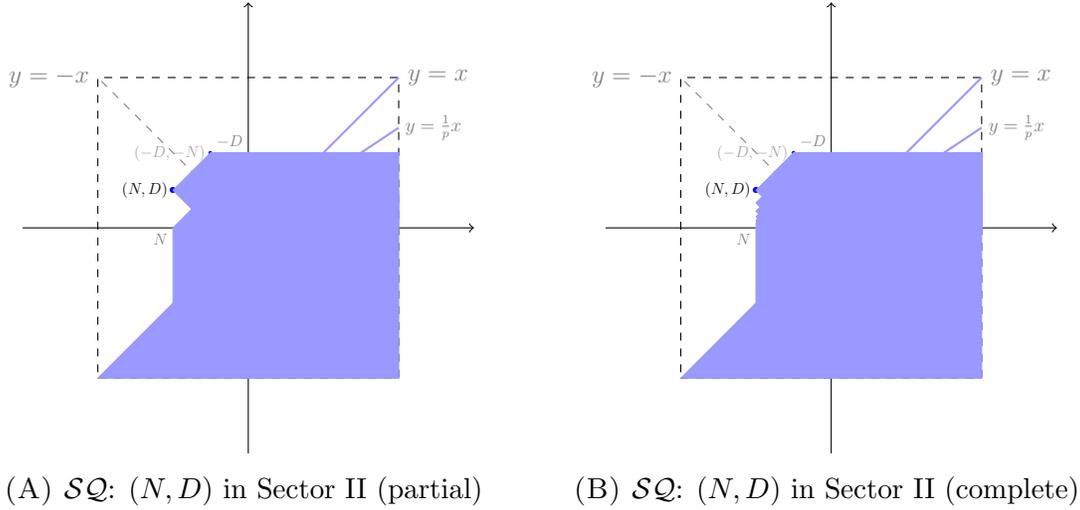

(A) $\mathcal{SQ}$: $(N, D)$ in Sector II (partial)

(B) $\mathcal{SQ}$: $(N, D)$ in Sector II (complete)

FIGURE 16. Work toward $\mathcal{SQ}$: $(N, D)$ in Sector II

However, to this we must add the vertical raindrops from $(N, D)$, each of which comes with its own version of the picture in Figure 16A; what results finally from the union as a safe subnormal quotient zone for $(N, D)$ in Sector II is as in Figure 16B.

9.8. $(N, D)$ **in Sector III.** To the safe quotient zone for $\mathcal{MID}$ quotients as in Figure 11B we may as usual add Sectors VI, VII, and VIII and the reflections of the special lines in Sector IV; it is easy to show that the horizontal line from $(N, D)$ extends to the right-hand boundary. When we include the vertical raindrops, those remaining in Sector III contribute their own horizontal line segments, while those in Sector II contribute their own version of Figure 16A. What results from the union is Figure 17 or as in Figure 4 on page 7.

## 10. CIRCUMSCRIBING SAFE ZONES

We do not claim to have obtained the maximal safe zones either for $\mathcal{MID}$ safe quotients or for subnormal safe quotients. A natural question, for example, is whether the quadrilaterals bounded by the horizontal and slant safe lines for a base point in Sector III can be added to the safe subnormal quotient zone (see Figure 17). This is not true, at least in general; with the aid of *Mathematica* [12] we can show that with $(N, D) = (-1/10, 2/5)$ the point $(2/5, 1/3)$ interior to one such quadrilateral is not a subnormal safe quotient. As well, there are other assorted routes to attempt to rule out points in the magic square from being safe, which we discuss and illustrate briefly.

First, a shift is $\mathcal{MID}$ if and only if its weights (or weights squared) are log completely alternating. If we fix a value of $p$, we may check using *Mathematica* [12] plots of the

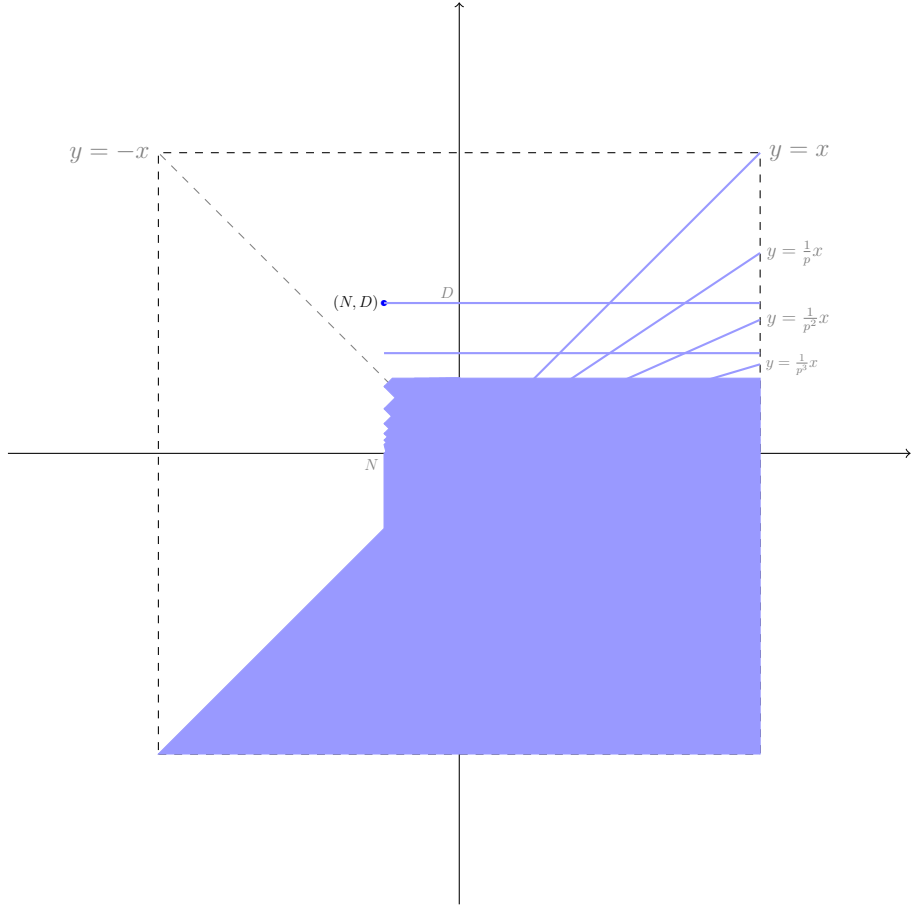

Figure 17. $\mathcal{SQ}$: $(N, D)$ in Sector III

areas of negativity for tests of the form

$$\sum_{i=0}^{n} (-1)^i \binom{n}{i} \ln\left(\alpha_{i+j}\right),$$

where the test starts at $\alpha_j$ and the test is for part of $n$-alternating for the log. The picture in Figure 18 gives (for $p = 3/2$, $(N, D) = (-1/2, -1/3)$) a version of such a plot. (A point is excluded if it is not on the same "side" of the curve as the origin, or, equivalently, if it is not on the same side of the curve as the line segment which is part of the line $N = D$.) The picture is indicative of the fact that since the main diagonal is safe, it seems numerically difficult to exclude areas. We note that this problem is even more severe in the case of subnormal safe zones, as a glance at Figure 17 on page 24 will make clear.

Second, we may prove that certain zones are excluded. In particular, if $(M, P) \gg (N, D)$ (respectively, $(M, P) \gg_s (N, D)$), except in the trivial case $(M, P) = (N, D)$ we have $(N, D) \not\gg (M, P)$ (respectively, $(N, D) \not\gg_s (M, P)$) because of the partial order.

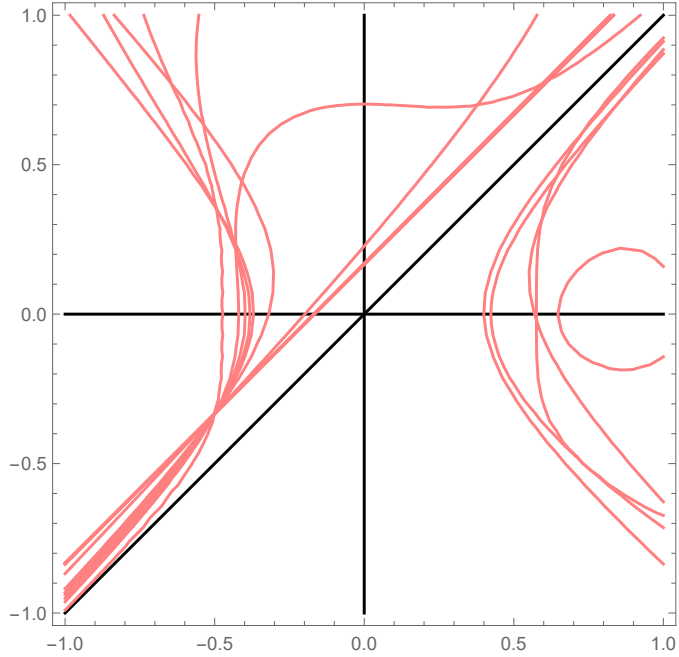

Figure 18. $\mathcal{MID}$ Exclusion using log completely alternating test curves

For example, and using transitivity, for $(N, D)$ in Sector I as in Figure 2A on page 6 we may exclude from its safe $\mathcal{MID}$ quotients every $(M, P)$ such that $P < D$ and $M$ puts the point above the diagonal line down and to the left from $(N, D)$. What results, in the case $p = 2$, $(N, D) = (-1/3, -1/2)$ is an excluded set as shown in Figure 19.

Third, we note that the excluded points form an open set (since either $\mathcal{MID}$ or subnormality, as the intersection of zones determined by weak inequality conditions, is a closed condition). So, for example, for $\mathcal{MID}$ safe quotients for a point $(N, D)$ in Sector II as in Figure 2B on page 6, one can show that the open segment between $(N, D)$ and $(N, -D)$ is not safe (the relevant quotients fall in sectors not yielding even subnormality). This means that for any compact subset of this open segment, there is actually an open "tube" of failure points.

In the paper [6], forthcoming, we will explore a notion which, applied in the particular case of the GRWS quotients, allows for further exclusion from quite a different point of view.

**Remark 10.1.** First, the safe zones we have found for $\mathcal{MID}$ quotients would appear to be independent of $p$ and depend only on the pair $(N, D)$; we conjecture that this is true for the full safe zones but do not understand why this should hold. It surely does not for subnormal safe quotients, because of the special lines which are clearly dependent upon $p$.

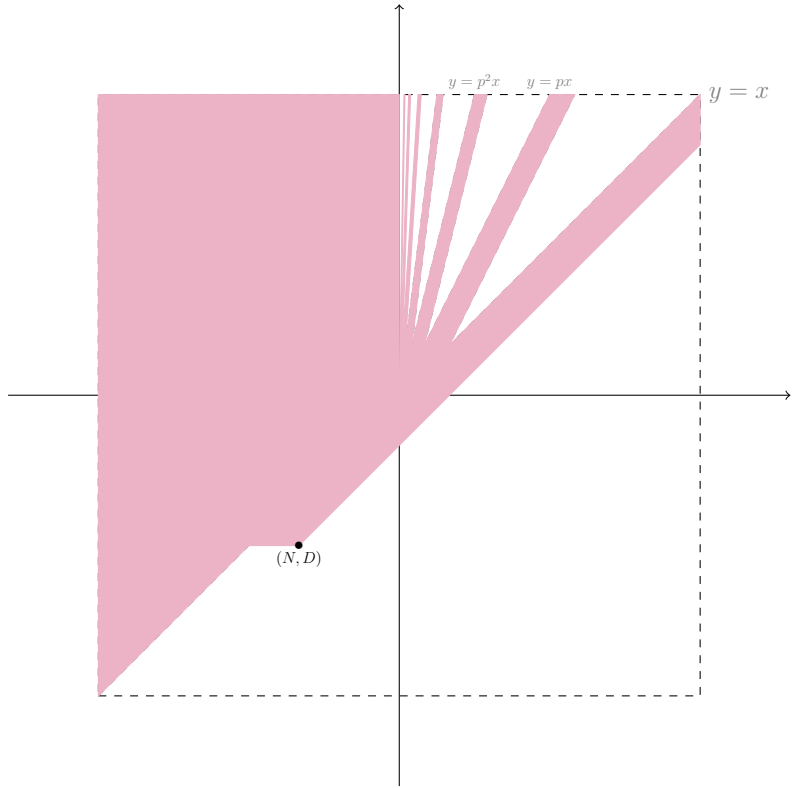

FIGURE 19. Excluded zone for a base point in Sector VIII

We observe next that the partial orders, here used within the class of GRWS, can in fact (with due care) be extended to the class of all shifts, and we will take this up in a subsequent paper [7].

**Acknowledgments.** The authors wish to express their appreciation for the support and warm hospitality during various visits (which materially aided this work) to Bucknell University, the University of Iowa, and the Université des Sciences et Technologies de Lille, and particularly to the Mathematics Departments of these institutions. The second named author was partially supported by NSF grant DMS-2247167. Several examples in this paper were obtained using calculations with the software tool *Mathematica* [12].

**Declarations**

<u>Conflict of Interest declaration</u>: The submitted work is original. It has not been published elsewhere in any form or language (partially or in full), and it is not under simultaneous consideration or in press by another journal. The authors have no competing interests to declare that are relevant to the content of this article.

<u>Data availability:</u>  The manuscript has no associated data.

## References

[1] J. Agler, Hypercontractions and subnormality, *J. Operator Theory* **13**(1985), 203–217.

[2] A. Athavale, On completely hyperexpansive operators, *Proc. Amer. Math. Soc.* **124**(1996), 3745–3752.

[3] C. Benhida, R.E. Curto, and G.R. Exner, Moment infinitely divisible weighted shifts, *Complex Analysis and Operator Theory* **13**(2019), 241–255.

[4] C. Benhida, R.E. Curto and G.R. Exner, Moment infinite divisibility of weighted shifts: Sequence conditions, *Complex Analysis and Operator Theory* **16: 1**(2022), art. 5, 1–23.

[5] C. Benhida, R.E. Curto and G.R. Exner, Geometrically Regular Weighted Shifts, preprint, arXiv:2309.0588.

[6] C. Benhida, R.E. Curto, and G.R. Exner, Signed representing measures (Berger charges) in subnormality and related properties of weighted shifts, forthcoming.

[7] C. Benhida, R.E. Curto, and G.R. Exner, Partial orders for Schur quotients of weighted shifts, forthcoming.

[8] C. Berg, J.P.R. Christensen and P. Ressel, *Harmonic Analysis on Semigroups*, Springer Verlag, Berlin, 1984.

[9] J. Bram, Subnormal operators, *Duke Math. J.* **22**(1965), 75–94.

[10] J.B. Conway, *The Theory of Subnormal Operators*, Mathematical Surveys and Monographs, vol. 36, Amer. Math. Soc., Providence, 1991.

[11] R. Gellar and L.J. Wallen, Subnormal weighted shifts and the Halmos-Bram criterion, *Proc. Japan Acad.* **46** (1970), 375–378.

[12] Wolfram Research, Inc. *Mathematica, Version 12.1*, Wolfram Research Inc., Champaign, IL, 2019.